\documentclass[12pt,thmsa]{article}
\usepackage{amsfonts}
\usepackage{amssymb}
\newtheorem{teo}{Theorem}[section]

\newtheorem{obs2}[teo]{Remark}

\newtheorem{tea}{Theorem}[subsection]

\newtheorem{no2}[teo]{Note}

\newtheorem{no3}[tea]{Note}

\newcommand{\Gal}{{\rm Gal}}

\newcommand{\Q}{\mathbb{Q}}

\begin{document}
\title{{\bf Improvements on Dieulefait-Manoharmayum and applications
}}
\author{Luis Dieulefait
\\
Dept. d'\'{A}lgebra i Geometria, Universitat de Barcelona;\\
Gran Via de les Corts Catalanes 585;
08007 - Barcelona; Spain.\\
e-mail: ldieulefait@ub.edu\\
 }
\date{\empty}

\maketitle

\vskip -20mm
%\titlerunning{Images of $3$-dimensional Galois representations
%

\begin{abstract}
We improve the results in a previous article of Dieulefait and Manoharmayum and we deduce some stronger modularity results.
\end{abstract}

\section{The results}

First let us prove the following lemma for $\rho$ an odd, irreducible, $2$-dimensional Galois representations of the absolute Galois group of $\Q$, with values in a finite field of odd characteristic $p$:\\

\bf{Lemma:} \rm Let $M$ be the quadratic number field ramifying at $p$ only. If $\rho$ is dihedral induced from $M$,  the Serre's weight of $\rho$ is $k < p$ and its Serre's level is arbitrary, then $\rho$ restricted to $I_p$ (the inertia group at $p$) is reducible  and the action of $I_p$ (after passing to tame inertia) is given by level $1$ fundamental characters. \\

Proof: If $p$ is different from $2k-3$ this is proved in: [DM] (Dieulefait, Manoharmayum: ``Modularity of rigid Calabi-Yau threefolds over Q", in ``Calabi-Yau Varieties and Mirror Symmetry", Fields Institute Communications Series, American Mathematical Society 38 (2003) 159-166) for $k=4$ but the proof works in general (see the ``last comment" below). Also, this is why R. Taylor takes $p$ different from $2k-3$ in his modularity lifting result in the crystalline non-ordinary case (in his article ``On the meromorphic continuation of degree two L-functions", preprint). 

So, for $p= 2k-3$ let us try to see that a similar proof applies. Since $k$ is 
even, this $p$ is congruent to $1$ mod $4$, thus $M$ is a real quadratic field. 
Since $\rho$ is an odd representation and $M$ is real, when we restrict to $M$ the 
image must be in a SPLIT Cartan subgroup of $GL_2(\mathbb{F})$, where $\mathbb{F}$ is the field of coefficients of
$\rho$. If the inertia group 
$I_p$ acts through level $2$ fundamental characters, the cyclic group generated by the values of the corresponding character  $\psi_2^{k-1} = 
\psi_2^{(p+1)/2}$, which, precisely because $p=2k-3$, has order equal to $2$ in the 
projective representation, is not contained in this Cartan (because in the opposite case $M/\Q$ would not ramify at $p$), and thus corresponds to the quadratic extension $M/\Q$. But level $2$ fundamental characters have nothing to do with $M/\Q$, because it is $\chi^{(p-1)/2}$, where $\chi$ is the cyclotomic 
character, the character that corresponds to $M/\Q$. Moreover, since in this case we have that after restricting $\rho$ to the absolute Galois group of $M$ the action of $I_p$ is given just by (the group of scalar matrices generated by) $\chi$ , which is a character that extends to the absolute Galois group of $\Q$, it is clear that this contradicts the fact that the action of  $\rho$ restricted to $I_p$ is given by the character $\psi_2^{(p+1)/2}$.\\     
We get a contradiction, so $\rho$ restricted to $I_p$ must be reducible and given by level $1$ fundamental characters. \\

Of course, this Lemma has consequences: first those already noted in 
[DM], that when you exclude this degenerate case then a combination of 
modularity lifting results \`{a} la Diamond-Taylor-Wiles (more precisely, a result of Diamond-Flach-Guo or a similar result of R. Taylor) with those of Skinner- 
Wiles are enough in the crystalline case (as proved in [DM]: recall that this 
also uses a result of Breuil to guarantee ordinarity in the case of level $1$ fundamental characters 
under the hypothesis of crystalline lifts and $k<p$), thus we have: \\

\bf{Corollary 1} \rm (combining modularity liftings theorems (cf. [DM])): If $\rho$ is 
reducible or modular, of Serre's weight $k<p$ and any level, any crystalline irreducible lift 
of Hodge-Tate weights $(0,k-1)$ is modular. \\

See [DM] for a proof. The idea is just the following: in the case of level $1$ fundamental characters by a result of Breuil the deformation we are considering is ordinary, and then the results of Skinner-Wiles apply. In the case of level $2$ fundamental characters, the technical condition needed to apply results \`{a} la Diamond-Taylor-Wiles is satisfied, as proved in the lemma above.\\

Another consequence is that the principle of ``switching the residual characteristic" (used in the articles that prove some cases of Serre's conjecture by Dieulefait and Khare-Wintenberger)  holds for any level and weight: \\

\bf{Corollary 2} \rm (combine corollary 1 with ``existence of families" and ``lowering the conductor" or ``existence of minimal lifts"): If, for $k, N, p$ fixed with $k<p$ we know that: Serre's conjecture is true in characteristic $p$, weight $k$ and ANY level $N'$ dividing $N$; then for any prime $q >k$ Serre's conjecture in characteristic $q$, weight $k$ and level $N$ is true. (of course, we are taking $p$ and $q$ not dividing $N$).\\

(For example, this principle is used in my preprint  ``The level 1 weight 2 case of Serre's conjecture", to reduce the proof to the case of characteristic $p=3$, where the conjecture was proved by Serre for this weight and level in 1973. See also the preprint of Khare-Wintenberger ``On Serre's reciprocity conjecture for $2$-dimensional mod $p$ representations of $\Gal(\bar{\Q}/\Q)$" for a similar strategy). 
The proof is the same that in the mentioned papers on Serre's conjecture (``existence of minimal lifts" is proved in these preprints and ``existence of compatible families" in my paper ``Existence of compatible families and new cases of the Fontaine-Mazur conjecture", J. Reine Angew. Math. 577 (2004) 147-151) : if you start in characteristic $q$, after taking a minimal lift and building a compatible family containing it, you look at a $p$-adic member of this family. It only remains to prove modularity of this $p$-adic member, but if you assume Serre's conjecture in characteristic $p$ (the weight is fixed at $k$ and the level bounded  by $N$ in all the process, just observe that the level may descend after switching, because sometimes conductors descend when reducing mod $p$) the modularity of this $p$-adic representation follows from corollary 1.\\
 
Finally, just observe that in the article [DM] to prove modularity of rigid 
Calabi-Yau threefolds we have not a very good result at $p=5$ precisely because 
for $k=4$, $5= 2k-3$ (this is the reason why in [DM] we obtained a better result at $p=7$). 
Now, with the above lema and corollary 1, we also get as in [DM] from the truth of Serre's conjecture on the field of $5$ elements: \\

\bf{Corollary 3:} \rm Any rigid Calabi-Yau threefold defined over $\Q$ with good reduction 
at $5$ is modular. \\

Last comment: in [DM] the full proofs of the above lemma, for $p$ different than $2k-3$, 
and of the above corollary 1, are given. They are given for $k=4$ and it is even 
explicitely said that proofs work for general $p$ and $k$ under the condition 
 $(p+1)/ \gcd(k-1,p+1) >2$  and $p > k$, which is the same as saying $p >k$ and 
different than $2k-3$.\\
 
\bf{Epilogue}: \rm the criterion ``Any rigid Calabi-Yau threefold over $\Q$ with good reduction at $3$ is modular" also follows, under the assumption that the modularity lifting result of Diamond-Flach-Guo (in their paper ``Adjoint motives of modular forms and the Tamagawa number conjecture", Ann Sci ENS. 37 (2004), no. 5, 663-727) can be extended, with the rest of conditions unchanged, to the case of crystalline representations of weights $(0, p)$ where $p$ (an odd prime) is the residual characteristic (currently the proof given by Diamond-Flach-Guo assumes that the weights $(0,w)$ satisfy $w < p$).\\
Remark: Recall that the $p$-adic Galois representation attached
 to a rigid Calabi-Yau threefold is crystalline of weights $(0,3)$ for any prime $p$ where the variety has good reduction.\\
  The proof is similar to the proof of corollary 3 above: as in Wiles original paper, we start by observing that the mod $3$ representation is either modular or reducible. In the ordinary case the results of Skinner and Wiles suffice for a proof. In the non-ordinary case, results of Berger-Li-Zhu (in  ``Construction of some families of 2-dimensional crystalline representations", Math. Annalen 329 (2004), no. 2, 365-377) give a precise description of the action of $I_p$ in the residual mod $3$ representation $\rho$: it acts through level $2$ fundamental characters and with Serre's weight equal to $2$. So again we can apply the lemma above, since $\rho$ is irreducible (fundamental characters of level $2$) and $k=2 < 3$, and conclude that $\rho$ restricted to the quadratic field ramifying only at $3$  is still irreducible. This is the technical condition required to apply the ``stronger version" of the result of Diamond-Flach-Guo in the case of a crystalline deformation of weights $(0,p)$.\\ 
(Hope: given the precise information about the action of $I_p$ on the residual representation in the case of weights $(0,p)$, non-ordinary, obtained by Berger-Li-Zhu, perhaps the extension of the result of Diamond-Flach-Guo to cover also this case is something that may be obtained in the near future).\\

\bf{Corollary 4:} \rm Assume that the result of Diamond-Flach-Guo can be extended to the case of crystalline deformations of weights $(0,p)$.  Then any rigid Calabi-Yau threefold defined over $\Q$ with good reduction at $3$ is modular.

\end{document}